\documentclass[10pt,leqno,twoside]{amsart}
\numberwithin{equation}{section}
\usepackage{latexsym}
\usepackage{amsmath}
\usepackage{amssymb}
\usepackage{mathrsfs}
\usepackage{graphicx}
\usepackage{ifthen}
\usepackage[T1]{fontenc}
\usepackage[latin1]{inputenc}

\numberwithin{equation}{section}

\usepackage{latexsym}
\usepackage{amsmath}
\usepackage{amssymb}
\newcommand{\cA}{{\mathcal A}}

\newcommand{\cD}{{\mathcal D}}

\newcommand{\cE}{{\mathcal E}}
\newcommand{\cB}{{\mathcal B}}
\newcommand{\cSM}{{\mathcal {SM} }}
\newcommand{\C}{{\mathbb C}}
\newcommand{\CC}{{\mathbb C}}
\newcommand{\EE}{{\mathbb E}}
\newcommand{\N}{{\mathbb N}}
\newcommand{\R}{{\mathbb R}}
\newcommand{\RR}{{\mathbb R}}
\newcommand{\PP}{{\mathbb P}}
\newcommand{\ve}{{\varepsilon}}

\newcommand{\nn}{\nonumber}

\begin{document}

\title[Dynamics of Nematic Liquid Crystal Flows]
{Thermodynamical Consistent Modeling and Analysis of  Nematic Liquid Crystal Flows}

\author{Matthias Hieber}
\address{Technische Universit\"at Darmstadt\\
        Fachbereich Mathematik\\
        Schlossgarten-Strasse 7\\
        D-64289 Darmstadt, Germany}
\email{hieber@mathematik.tu-darmstadt.de}

\author{Jan Pr\"uss}
\address{Martin-Luther-Universit\"at Halle-Witten\-berg\\
         Institut f\"ur Mathematik \\
         Theodor-Lieser-Strasse 5\\
         D-06120 Halle, Germany}
\email{jan.pruess@mathematik.uni-halle.de}


\begin{abstract}
The general Ericksen-Leslie model for the flow of nematic liquid crystals is reconsidered in the non-isothermal case aiming for thermodynamically consistent models.
The non-isothermal simplified model is then investigated analytically. A fairly complete dynamic theory is developed by
analyzing these systems as  quasilinear parabolic evolution equations in an $L_p-L_q$-setting.
First, the existence of a unique, local strong solution is proved. It is then shown that this solution extends to a global strong solution provided the initial data are close to an
equilibrium or the solution is eventually bounded in the natural norm of the underlying state space. In these cases the solution converges exponentially to an equilibrium in the
natural state manifold.
\end{abstract}

\maketitle

\section{Introduction}

The continuum theory of liquid crystals was developed by Ericksen and Leslie during the 1960's in  their pioneering work \cite{Eri62,Les68}. This theory models nematic liquid crystal flow
from a hydrodynamical point of view  and reduces to the Oseen-Frank theory in the static case, see \cite{Ose33},\cite{Fra58}. It describes the evolution of the complete system
under the influence of the velocity $u$ of the fluid and the orientation configuration $d$ of rod-like liquid crystals. Hence, $d=d(t,x)$ is a unit vector in $\R^3$.
The original derivation \cite{Eri62,Les68} is based on the conservation  laws for mass and linear as well as angular momentums.
General liquid crystal materials  are described by the Landau-de Gennes theory \cite{DG74} from a unified point of view.

The Ericksen-Leslie theory is nowadays widely used as a model for the flow of nematic liquid crystals, see for example the works of
Ericksen and Kinderlehrer \cite{EK87},  Chandrasekhar, \cite{Cha92}, deGennes and Prost \cite{DGP95} as well as Virga \cite{Vir94}.

Note that these models are mostly formulated in an isothermal environment and are, in general, neither thermodynamically consistent nor
thermodynamically stable. To the best of our our knowledge, only very few articles are dealing so far  with the {\em thermodynamical consistency} of these models.
Concerning the Ericksen-Leslie model, for a physically rigorous derivation we refer to the work of M\"uller \cite{Mue85} concentrating on the modeling aspect,
and for recent analytical work to Feireisl, Rocca, Schiperna \cite{FRS11}, Feireisl, Fr{\'e}mont, Rocca, Schiperna \cite{FFRS12} and Li, Xin \cite{LX14}. Non-isothermal Landau-De Gennes nematic liquid crystal
flows were investigated in the recent articles \cite{FSRZ12} and \cite{FSRZ13}.

The aim of this  paper is twofold: first, we reconsider the Ericksen-Leslie approach from the perspective of  thermodynamical consistency and stability.
Following arguments from thermodynamics and employing entropy principles, we derive consistent models in a mathematically efficient way, even in the case of
compressible fluids. Let us emphasize that, in the end, our model contains the classical Ericksen-Leslie model in its general form as a special case.

Secondly, we investigate our model analytically.  Restricting ourselves to the case of constant density and not taking into account so called stretching, we develop a rather
complete dynamic theory for the equations representing these models. More precisely, we first prove  the existence of a unique, local strong solution to this system.
We further show that this solution extends to a global, strong solution, provided the initial data are close to an equilibrium or the solution is eventually bounded in the
natural norm of the underlying state space. In this case the solution converges exponentially to an equilibrium in the natural state manifold.
The results obtained thus parallel those proved recently by Hieber, Nesensohn, Pr\"uss and Schade in \cite{HNPS14} dealing with the isothermal situation.
For results concerning the asymptotic behaviour of solutions
in the situation of the whole space $\R^3$,  we refer to the work of Dai and Schonbek \cite{DS14}.

The nowadays called simplified Ericksen-Leslie model in the isothermal situation was introduced and investigated  first by Lin in \cite{Lin89} and \cite{Lin91}. Lin and Liu
\cite{LL95}, \cite{LL00}
studied the situation,  where the nonlinearity in the equation for the director $d$ is replaced by a Ginzburg-Landau energy functional. The existence of global weak solutions
to this system  in dimension $2$ or $3$ was proved under suitable assumptions on the intial data.  For related results see \cite{HW10} and \cite{LLW10}.
Wang proved in \cite{Wan11} global well-posedness for the simplified  system for initial data being small in $BMO^{-1}\times BMO$
in the case of a whole space by combining techniques of Koch and Tataru with methods from harmonic maps to certain Riemannian manifolds.

The general Ericksen-Leslie model (in the isothermal situation) is based on  the  Oseen-Frank energy density functional which takes
into account stretching as well as rotational effects for the director field. In the special case of homogeneous isotropic elasticity the equation for $d$ reads  as
\begin{equation}\label{eqmh:d}
\partial_td  + u \cdot \nabla d  - Vd + \frac{\lambda_2}{\lambda_1}Dd = -\frac{1}{\lambda_1}(\Delta d + |\nabla d_2^2d) + \frac{\lambda_2}{\lambda_1}(Dd \cdot d)d
\quad \text{in } (0,T) \times \Omega.
\end{equation}
Here $D=\frac12(\nabla u+[\nabla u]^T)$ denotes the symmetric, $V=\frac12(\nabla u- [\nabla u]^{\sf T})$ the anti-symmetric part of the deformation tensor and
$\lambda_1,\lambda_2 \in \R\backslash\{0\}$ are material coefficients. Modifications of this model were considered by Coutard and Shkoller
in \cite{CS01} in which the above equation \eqref{eqmh:d} for d  is replaced by a Ginzburg-Landau type approximation.
\begin{equation}\label{eq:stretch}
 \gamma(\partial_td  + u \cdot \nabla d  + d \cdot \nabla u)= \Delta d - \frac{1}{\varepsilon^2}(|d|_2^2 -1)d \quad  \text{in } (0,T) \times \Omega.
\end{equation}
They proved local wellposedness for this system  as well as a  global existence result for small data within this setting.
Note, however, that in this case the presence of the stretching term $d \cdot \nabla u$ causes loss of total energy balance and, moreover, the
condition $|d|_2= 1$  in  $(0,T) \times \Omega$, is not preserved anymore.
For recent results on the general Ericksen-Leslie model with $d$ satisfying \eqref{eqmh:d}, we refer to the articles
\cite{LLW10}, \cite{HLW13}, \cite{WZZ13}, \cite{LWX13}, \cite{MGL14}, \cite{HLX14} and \cite{FL15},  which contain well-posedness criteria  for the
general system under various assumptions on the Leslie coefficients.
For results on non-isothermal models including the above mentioned stretching term,  see \cite{FFRS12} and references therein.

Let us stress at this point that an important novelty of our approach lies in the fact that the complete model described in Section 2.7 is rigoroulsy proven to be thermodynamically
consistent and stable. Specialising to the isothermal situation, we rediscover in particular the classical general Ericksen-Leslie system.
It is interesting to compare our approach with the approach of M\"uller \cite{Mue85}, and with the
energy variational approach developed by  Liu and coworkers \cite{LWX13} and by Virga \cite{Vir94}.

The plan for this contribution is as follows: Section 2 is devoted to the modeling of liquid crystals.
In particular, based on the entropy principle, we derive a model of Ericksen-Leslie type which is thermodynamically consistent and stable.
In Section 3 the equilibria of the system are identified - which are zero velocities and constant temperature and director - and it is proved that these are thermodynamically stable.
The negative total entropy is shown to be strict Ljapunov functional, in particular the model is thermodynamically consistent, In Section 4 we prove local well-posedness of the non-isothermal simplified model and construct the resulting local semiflow in the natural state manifold of the system.
We show that each solution which does not
develop singularities in a sense to be specified converges to a unique equilibrium. These results are proved by means of techniques involving maximal $L_p$-regularity
and results on quasilinear paprabolic evolution equations.  For these methods, we refer to the booklet  by Denk, Hieber, and Pr\"{u}ss \cite{DHP03} and to the work of
Pr\"{u}ss and Simonett \cite{PS04}, K\"{o}hne, Pr\"{u}ss and Wilke \cite{KPW10},  and  LeCrone, Pr\"{u}ss and Wilke \cite{LPW14}.

By means of these techniques we are also able to prove analogous results for the full model, which, however, due to limitation of space will be presented  elsewhere.

\section{Thermodynamical Consistent Modeling }
\noindent
In this section we aim to give a self-contained presentation of a {\em thermodynamically consistent modeling of liquid crystals}.
Like this we are able to refrain from  refering to the orginal papers \cite{Eri62}, \cite{Les68}, which are not easily accessible to a mathematical audience.
Our approach does not only  extend the classical Ericksen-Leslie model to the non-isothermal situation in a thermodynamical consistent and stable  way but it also allows
to exhibit the  physical and mathematical beauty of this model.

In this section, $\Omega\subset\RR^n$ always denotes a domain with $C^1$-boundary.

\bigskip

\noindent
{\bf 1. First Principles}\\
We begin with the balance laws  of mass, momentum, and energy. They read as
\begin{align}
\partial_t \rho +{\rm div}(\rho u)&=0\quad &\mbox{in } \Omega,\nn\\
\rho(\partial_t +u\cdot\nabla)u +\nabla \pi &= {\rm div}\, S \quad &\mbox{in } \Omega,\\
\rho(\partial_t +u\cdot\nabla)\epsilon +{\rm div}\, q &= S:\nabla u -\pi{\rm div}\, u\quad &\mbox{in } \Omega,\nn\\
u=0,\quad q\cdot \nu&=0 \quad &\mbox{on } \partial\Omega.\nn
\end{align}
Here $\rho$ means density, $u$ velocity, $\pi$ pressure, $\epsilon$ internal energy, $S$ extra stress and $q$ heat flux. This immediately gives conservation of the total energy.
In fact, we have
$$\rho(\partial_t+u\cdot\nabla) e +{\rm div}(q+\pi u-Su)=0 \quad \mbox{ in } \Omega,$$
where $e:= |u|^2/2 +\epsilon$ means the total mass specific energy density (kinetic  and internal). The energy flux $\Phi_e$ is given by $\Phi_e:=q+\pi u -Su$. Integrating over $\Omega$ yields
$$\partial_t {\sf E}(t) =0,\quad  {\sf E}(t) = {\sf E}_{kin}(t) + {\sf E}_{int}(t) =\int_\Omega \rho(t,x) e(t,x)dx,$$
provided
\begin{equation}\label{bc}
q\cdot\nu=u=0 \quad \mbox{ on }\partial\Omega.
\end{equation}
Hence, if \eqref{bc} holds, total energy is preserved, independent of the particular choice of $S$ and $q$.

\bigskip

\noindent
{\bf 2. Thermodynamics}\\
Assume a given free energy $\psi$  of the form $\psi= \psi(\rho,\theta,\tau)$, where $\theta$ denotes the (absolute) temperature and $\tau$ will be specified  later.
We then have the following thermodynamical relations:
\begin{align}
\epsilon &= \psi +\theta \eta \quad \mbox{ internal energy},\nn\\
\eta &= -\partial_\theta \psi \quad \mbox{ entropy},\\
\kappa &= \partial_\theta \epsilon = -\theta\partial_\theta^2\psi \quad \mbox{ heat capacity} \nn
\end{align}
Later on, for well-posedness of the heat problem, we require $\kappa>0$, i.e.\ $\psi$ to be  strictly concave with respect to  $\theta\in(0,\infty)$.

In the classical case, where $\psi$ depends only on $\rho$ and $\theta$, we have the {\em Clausius-Duhem equation}
$$\rho (\partial_t+u\cdot\nabla) \eta +{\rm div}(q/\theta) = S:\nabla u/\theta -q\cdot\nabla \theta/\theta^2+(\rho^2\partial_\rho\psi-\pi)({\rm div}\, u)/\theta\quad \mbox{in }\Omega.$$
Hence, in this case  the entropy flux $\Phi_\eta$ is given by $\Phi_\eta:=q/\theta$ and the entropy production by
$$ \theta r:=S:\nabla u -q\cdot\nabla \theta/\theta+(\rho^2\partial_\rho\psi-\pi)({\rm div}\, u).$$
Employing  the boundary conditions \eqref{bc}, we obtain for the total entropy ${\sf N}$ by integration over $\Omega$
$$ \partial_t {\sf N}(t) = \int_\Omega r(t,x) dx\geq0,\quad {\sf N}(t)= \int_\Omega \rho(t,x)\eta(t,x) dx,$$
provided $r\geq 0$ in $\Omega$. As ${\rm div}\, u$ has no sign we require
\begin{equation}\label{Maxwell}
 \pi=\rho^2\partial_\rho\psi,
 \end{equation}which is the famous {\em Maxwell relation}.
Further, as $S$ and $q$ are independent, this requirement leads to the classical conditions
\begin{equation}\label{EP}
S:\nabla u\geq 0 \quad \mbox{ and }\quad q\cdot\nabla \theta\leq 0.
\end{equation}
Summarizing, we see that whatever one chooses for $S$ and $q$, one always has conservation of energy and the total entropy is non-decreasing provided \eqref{EP},\eqref{Maxwell} and \eqref{bc}
are satisfied. Thus, these conditions ensure the thermodynamic consistency of the model.

As an example for $S$ and $q$ consider the classical laws due to Newton and Fourier which are given by
$$ S:=S_N:= 2\mu_s D +\mu_b {\rm div}\, u\,I, \quad 2D=(\nabla u +[\nabla u]^{\sf T}), \quad q = -\alpha_0 \nabla \theta.$$
In this case, \eqref{EP} is satisfied as soon as $\mu_s\geq0$, $2\mu_s+n\mu_b\geq0$ and $\alpha_0\geq0$ hold. Note that it does not matter at all
whether $\mu_s,\mu_b, \alpha_0$ are constants or whether they depend on $\rho,\theta$, or on other variables.

\bigskip

\noindent
{\bf 3. Nematic Liquid Crystals}\\
For isotropic nematic liquid crystals we assume a free energy density $\psi$ of the form
$$ \psi=\psi(\rho,\theta,\tau),\quad \mbox{ with } \tau = \frac{1}{2}|\nabla d|_2^2.$$
Here $d$ means the orientation vector, also called the {\em director}, which should satisfy the condition
$$
|d|_2^2 := \sum_{j=1}^n d_j^2=1.
$$
Note that $2\tau ={\rm tr}(\nabla d[\nabla d]^{\sf T})$ is the first invariant of the matrix $\nabla d[\nabla d]^{\sf T}$, and for its last invariant it holds
${\rm det}(\nabla d[\nabla d]^{\sf T}) =({\rm det }\nabla d)^2=0$, as $\nabla d d =0$ by $|d|_2=1$.

We neglect spin energy below but take into account transport of energy due to couple stress. This means that the energy flux is replaced by
$$\Phi_e:=q+\pi u -Su- \Pi\cD_td, \quad \cD_t=\partial_t+u\cdot\nabla d,$$
where $\Pi$ has to be modeled.

As constitutive laws we will employ
\begin{align}
S= S_N + S_E + S_L, \quad S_E= -\lambda \nabla d[\nabla d]^{\sf T},\quad q= -\alpha_0\nabla\theta -\alpha_1(d\cdot\nabla\theta) d.
\end{align}
$S_N$ means the {\em Newton stress} introduced above, $S_E$ the {\em Ericksen stress}, and $S_L$ the {\em Leslie stress} which will be defined later.
Assuming these two constitutive  laws we derive in the following the balance of entropy, i.e.\ the {\em Clausius-Duhem equation}.
A short computation gives
\begin{align}
\rho(\partial_t+u\cdot\nabla) \eta + {\rm div}\, \Phi_\eta = r,
\end{align}
with $\Phi_\eta = q/\theta$, and
\begin{align*}
\theta r &= -q\cdot\nabla\theta/\theta + 2\mu_s|D|_2^2 + \mu_b|{\rm div}\, u|^2 + (\rho^2\partial_\rho\psi-\pi){\rm div}\, u\nn\\
&+ (\rho\partial_\tau\psi-\lambda) \nabla d[\nabla d]^{\sf T}:\nabla u + (\Pi-\rho\partial_\tau \psi\nabla d):\nabla\cD_td\\
&+ S_L:\nabla u +({\rm div}\Pi+\beta d)\cdot\cD_td.\nn
\end{align*}
for some scalar function $\beta$. Note that $d\cdot \cD_td=0$ as $|d|_2=1$, hence $\beta\in \RR$ can be chosen arbitrarily.

For the entropy production $r$ to be nonnegative, we require
$$  \mu_s\geq0,\quad 2\mu_s+n\mu_b\geq0,\quad \alpha_0\geq0,\quad \alpha_0+\alpha_1\geq0.$$
Except for the last one, these conditions are the well-known conditions from fluid dynamics, see Section 2.2. The subsequent terms in the definition of $r$ have no sign, hence we require them to vanish, which yields the relations
\begin{equation}\label{maxwel}
 \pi =\rho^2 \partial_\rho\psi,\quad \lambda = \rho\partial_\tau\psi\quad \Pi = \rho\partial_{\tau}\psi\nabla d.
\end{equation}
 Finally, to obtain nonnegativity of the last two terms, in the simplest case, we may assume that the Leslie stress $S_L$ vanishes, and
$$ \gamma \cD_td =  {\rm div}[ (\rho\partial_\tau\psi)\nabla]d+\beta d,$$
for some $\gamma=\gamma(\rho,\theta,\tau)\geq0$. The condition $|d|_2=1$ then requires $\beta= \lambda|\nabla d|^2,$
which leads to the equation
\begin{equation}\label{d-eq}
\gamma (\partial_t+u\cdot\nabla)d = {\rm div}[ \lambda\nabla]d+\lambda|\nabla d|^2d,
\end{equation}
a nonlinear convection-diffusion equation for $d$. This is the basic equation governing the evolution of the director field $d$. With these assumptions the entropy production reads as
$$\theta r = -q\cdot\nabla\theta/\theta + 2\mu_s|D|_2^2 + \mu_b|{\rm div}\, u|^2 +  \frac{1}{\gamma}|{\sf a}|_2^2,$$
where
$${\sf a} = {\rm div}[ \lambda\nabla]d+\lambda|\nabla d|^2d = \gamma \cD_t d.$$
At the boundary $\partial\Omega$, energy should be preserved, which means $\Phi_e\cdot\nu=0$. As $q\cdot\nu=0$ and $u=0$ this yields
$$ \lambda \partial_\nu d \cdot \partial_t d=0.$$
This is clearly valid if $d$ satisfies the Neumann condition $\partial_\nu d=0$, which is physically reasonable.

\bigskip

\noindent
{\bf 4. Stretching and Vorticity}\\
Observe that the equation \eqref{d-eq} for $d$ admits the solutions $d=const$, no matter how the velocity field and the temperature field are defined. In this case the director
field is not at all affected by the fluid dynamics. This seems to be physically unrealistic and so the model should be adapted.

This can be done by introducing a so-called {\em stretching} stress. To introduce this stress we follow Leslie. Define
$P_d=I-d\otimes d$ the orthogonal projection onto $E_d:=\{d\}^\perp$, the vorticity $V$ according to $2V=\nabla u-[\nabla u]^{\sf T}$, and set
$$ {\sf n}= \mu_V Vd +\mu_D P_d Dd -\gamma{\cD}_t d,$$
where $\mu_V,\mu_D,\gamma$ are scalar functions of $\rho,\theta,\tau$ and  $\gamma>0$. For brevity we use the notation
$$  {\sf a}:= P_d {\rm div}(\lambda\nabla)d = {\rm div}(\lambda\nabla)d +\lambda|\nabla d|^2d.$$
Now we define the stretch tensor
\begin{equation}\label{Strech}
S_L^{stretch} = \frac{\mu_D+\mu_V}{2\gamma} {\sf n} \otimes d + \frac{\mu_D-\mu_V}{2\gamma} d \otimes {\sf n}.
\end{equation}
 This modification of the model does not change the entropy flux $\Phi_\eta=q/\theta$,
and the relevant entropy production becomes,
\begin{align*} S_L^{stretch}:\nabla u+ \cD_td\cdot{\sf  a}&= \frac{1}{\gamma}( |{\sf n}|^2 +\gamma\cD_t d\cdot {\sf n}))+{\sf a}\cdot\cD_td\\
& = \frac{1}{\gamma} ( |{\sf a}|^2 +({\sf n}+{\sf a})\cdot(\mu_V Vd +\mu_D P_dDd-{\sf a}).
\end{align*}
If we want to keep the total entropy production at the same level as in the previous section, the simplest way to achieve this is to set ${\sf n}+{\sf a}=0$, which yields the equation
\begin{equation}\label{stretched-eq}
\gamma(\partial_t d +u\cdot\nabla d) = {\rm div}(\lambda\nabla)d +\lambda|\nabla d|^2d + \mu_V Vd +\mu_D P_dDd.
\end{equation}
This is the stretched equation for $d$. Note that it preserves the constraint $|d|_2=1$. The entropy production is the same as before, we have
$$\theta r = [\alpha_0|\nabla\theta|_2^2+ \alpha_1(d|\nabla \theta)^2]/\theta + 2\mu_s|D|_2^2 + \mu_b|{\rm div}\, u|^2 +\frac{1}{\gamma}|P_d {\rm div}(\lambda\nabla)d|_2^2.$$
In particular, ${\sf N}$
satisfies
$\partial_t {\sf N}(t) = \int_\Omega r(t,x)dx,$
and so $-{\sf N}$ will shown below to be a strict Lyapunov functional for the system, as soon as
\begin{equation}
 \mu_s>0,\quad 2\mu_s+n\mu_b>0,\quad \alpha_0>0,\quad \alpha_0+\alpha_1>0, \quad \gamma>0,
 \end{equation}
and
\begin{equation}
 \kappa>0,\quad \lambda>0,\quad \partial_\rho\pi>0.
\end{equation}
Note that no conditions on the new parameter functions $\mu_D,\mu_V$ are needed, so far.

\bigskip

\noindent
{\bf 5. Additional Dissipation}\\
We may add additional dissipative terms in the stress tensor of the form
\begin{equation}\label{Dissipation}
 S_L^{diss} =\frac{\mu_p}{\gamma} ({\sf n}\otimes d + d\otimes {\sf n}) + \frac{\gamma\mu_L+\mu_P^2}{2\gamma}(P_dDd\otimes d + d\otimes P_dDd) +\mu_0 (Dd|d)d\otimes d,
\end{equation}
where as before $2D= \nabla u+[\nabla u]^{\sf T}$ and $2V= \nabla u-[\nabla u]^{\sf T}$ are the symmetric and antisymmetric parts of the rate of strain tensor $\nabla u$.
Note that the tensor $S_L^{diss}$ is symmetric.
Adding these terms to the stress tensor will be thermodynamically consistent provided their contribution to the entropy production ensures that the total entropy production remains nonnegative. By a simple calculation we obtain
\begin{equation}\label{Diss-Entr}
 S_L^{diss}:\nabla u = S_L^{diss}:D = 2\frac{\mu_P}{\gamma}({\sf n}|P_dDd) + (\mu_L+\frac{\mu_P^2}{\gamma})|P_dDd|_2^2 + \mu_0 (Dd|d)^2.
\end{equation}
So with ${\sf n} =-{\sf a}$, the total relevant dissipation amounts to
\begin{align*}
&\frac{1}{\gamma}(|{\sf a}|_2^2 + 2\mu_P({\sf n}|P_dDd) + \mu_P^2|P_dDd|_2^2) + \mu_L |P_d Dd|_2^2 +\mu_0 (Dd|d)^2\\
&= \frac{1}{\gamma}|{\sf a}-\mu_PP_dDd|_2^2 + \mu_L |P_d Dd|_2^2 +\mu_0 (Dd|d)^2,
\end{align*}
hence the total entropy production becomes
\begin{align*} \theta r &=[\alpha_0|\nabla\theta|_2^2+ \alpha_1(d|\nabla \theta)^2]/\theta + 2\mu_s|D|_2^2 + \mu_b|{\rm div}\, u|^2 \\
&+\frac{1}{\gamma}|P_d{\rm div}(\lambda\nabla)d-\mu_P P_dDd|_2^2+\mu_L |P_d Dd|_2^2 +\mu_0 (Dd|d)^2.
\end{align*}
Note that so far the parameter functions $\mu_j$, $j=0,s,b,V,D,P,L$, $\alpha_0,\alpha_1$, and $\gamma$ for thermodynamical consistency are only subject to the requirements
\begin{align}\label{thermod-consist}
\alpha_0,\alpha_0+\alpha_1\geq0, \quad \mu_s,  2\mu_s+n\mu_b\geq0,\quad \mu_0,\mu_L\geq0,\quad \gamma>0.
\end{align}
Recall that all parameters functions are allowed to be functions of $\rho,\theta,\tau$.

\bigskip

\noindent
{\bf Remark.} {\bf (i)} A more refined algebra shows that it is enough to require
$$ 2\mu_s +\mu_L\geq0,\quad 2\mu_s +\mu_0\geq0$$
in the incompressible case, and additionally
$$ \frac{\mu_0^2}{n^2}\leq (2\mu_s + \mu_0)( \frac{2\mu_s}{n} +\mu_b + \frac{\mu_0}{n^2})$$
in the compressible case.

\medskip

\noindent
{\bf (ii)} We  want to stress that in case $\mu_V=\gamma$, our parameters $\mu_s,\mu_0,\mu_V,\mu_D,\mu_P,\mu_L$ are in one-to-one correspondence to the famous Leslie parameters $\alpha_1,\ldots,\alpha_6$. This shows that our model contains the isotropic Ericksen-Leslie model as a special case.

\bigskip

\noindent
{\bf 6. Conservation of Angular Momentum in 3D}\\
We briefly discuss conservation of momentum in the physical important three-dimensional case. Recall that the mass specific density ${\sf m}$ of angular momentum is defined by
$$ {\sf m} = x \times u.$$
Balance of angular momentum reads as follows
$$ \cD_t (\rho{\sf m}) +{\rm div}( \rho u\times {\sf m} - x\times T) = -e_i\times T_i,$$
where we use Einstein's sum convention and $T_i$ denotes the $i$-th row of the stress tensor $T$.
Thus the flux of angular momentum $\Phi_{m}$ is given by
$$ \Phi_{m} = \rho u\times {\sf m} - x\times T.$$
It is well-known that $e_i\times T_i=0$ in case $T$ is symmetric. Therefore we may concentrate on the non-symmetric part of $T$ which is given by
$$T^{as}= \frac{\mu_V}{2\gamma} ( {\sf n}\otimes d-d\otimes {\sf n}).$$
This implies
\begin{align*}
e_i\times T^{as}_i &= \frac{\mu_V}{2\gamma}( {\sf n}\times d -d\times{\sf n}) =  -\frac{\mu_V}{\gamma}d\times{\sf n}\\
&=  \frac{\mu_V}{\gamma} d\times {\rm div}(\lambda\nabla) d\\
&= \partial_i(\frac{\mu_V}{\gamma}\lambda d\times \partial_id) - \partial_i(\frac{\mu_V}{\gamma}) \lambda d\times\partial_id - \frac{\mu_V\lambda}{\gamma} \partial_id\times\partial_i d\\
&=  \partial_i(\frac{\mu_V}{\gamma}\lambda d\times \partial_id) - \partial_i(\frac{\mu_V}{\gamma}) \lambda d\times\partial_id.
\end{align*}
This shows that $e_i\times T^{as}_i$ is a divergence provided $\nabla_x\frac{\mu_V}{\gamma}=0$, i.e.\ if
$$ \mu_V = c_0 \gamma,\quad \mbox{ for some constant } c_0\in\RR.$$
Then the flux of angular momentum becomes
$$ \Phi_{m} = \rho u\times {\sf m} - x\times T+ c_0\lambda d\times \nabla d.$$
We mention that if we also require so-called {\em objectivity} of the model, then $c_0=1$, which means $\mu_V=\gamma$.

\bigskip

\noindent
{\bf 7. The Complete Model: non-isothermal, compressible fluid, isotropic elasticity}\\
Summarizing, the complete model may be represented as
\begin{align}\label{cp-eqns}
\partial_t \rho +{\rm div}(\rho u)&=0\quad &\mbox{in } \Omega,\nn\\
\rho(\partial_t +u\cdot\nabla)u +\nabla \pi &= {\rm div}\, S \quad &\mbox{in } \Omega,\nn\\
\rho(\partial_t +u\cdot\nabla)\epsilon +{\rm div}\, q &= S:\nabla u -\pi{\rm div}\, u+ {\rm div}(\lambda\nabla d\cD_t d)\quad &\mbox{in } \Omega,\\
\gamma (\partial_t+u\cdot\nabla)d -\mu_V Vd- {\rm div}[ \lambda\nabla]d&=\lambda|\nabla d|^2d  + \mu_D P_d Dd,\quad &\mbox{in } \Omega,\nn\\
u=0,\quad q\cdot \nu=0,\quad \partial_\nu d&=0 \quad &\mbox{on } \partial\Omega.\nn\\
\rho(0)=\rho_0,\quad u(0)=u_0,\quad \theta(0)=\theta_0,\quad d(0)&=d_0\quad &\mbox{in } \Omega\nn
\end{align}
These equations have to be supplemented by the thermodynamical laws
\begin{align}\label{cp-thdyn}
&\epsilon = \psi +\theta \eta, \quad \eta = -\partial_\theta \psi, \quad \kappa = \partial_\theta \epsilon,\nn\\
& \pi =\rho^2 \partial_\rho\psi,\quad \lambda = \rho\partial_\tau\psi/\theta,
\end{align}
and by the constitutive laws
\begin{align}\label{cp-const}
S &= S_N + S_E + S_L^{stretch} + S_L^{diss}, \nn\\
S_N&= 2\mu_s D + \mu_b {\rm div}\, u \, I, \quad S_E=  -\theta\lambda \nabla d[\nabla d]^{\sf T},\\
S_L^{stretch}& = \frac{\mu_D+\mu_V}{2\gamma}{\sf n} \otimes d + \frac{\mu_D-\mu_V}{2\gamma} d \otimes {\sf n},\quad {\sf n}= \mu_V Vd +\mu_D P_dDd-\gamma\cD_td,\nn\\
S_L^{diss}& =\frac{\mu_P}{\gamma} ({\sf n}\otimes d + d\otimes {\sf n}) + \frac{\gamma\mu_L+\mu_P^2}{2\gamma}(P_dDd\otimes d + d\otimes P_dDd) +\mu_0 (Dd|d)d\otimes d,\nn\\
q&= -\alpha_0 \nabla\theta -\alpha_1(d|\nabla\theta)d.\nn
\end{align}
Here all coefficients $\mu_j,\alpha_j$ and $\gamma$ are functions of $\rho,\theta,\tau$. For thermodynamic consistency we require
\begin{align}\label{cp-consist}
\mu_s\geq0,\quad 2\mu_s+n\mu_b\geq0,\quad \alpha_0\geq0,\quad \alpha_0+\alpha_1\geq0,\quad \mu_0,\mu_L\geq0,\quad \gamma>0,
\end{align}
Finally,  we will use in addition the following conditions
\begin{align}\label{cp-stable}
&\mu_s>0,\quad 2\mu_s+n\mu_b>0,\quad \alpha_0>0,\quad \alpha_0+\alpha_1>0,\quad \gamma>0,\nn\\
& \kappa>0,\quad \lambda>0, \quad \partial_\rho \pi >0,
\end{align}
to identify the equilibria and to investigate their thermodynamic stability in Section 3.

\bigskip

\bigskip

\noindent
{\bf 8. The Complete Model: non-isothermal, compressible fluid, non-isotropic elasticity}\\
For the sake of completeness, we comment briefly on the non-isotropic case. Then $\psi=\psi(\rho,\theta,d,\nabla d)$, and the Ericksen stress tensor
 becomes $S_E = -\rho \frac{\partial\psi}{\partial \nabla d}[\nabla d]^{\sf T}$. Following the derivation in Sections 2.3 and 2.4, here the energy and entropy fluxes read again as
$$\Phi_e:=q+\pi u -Su- \Pi\cD_t d,\quad \Phi_\eta = q/\theta,$$
and the equation for $d$ becomes
$$ \gamma\cD_t d = P_d {\sf a} ,\quad {\sf a} = \partial_i(\rho\nabla_{\partial_id}\psi)-\rho \nabla_d \psi,$$
in the case without stretching, and
$$ \gamma\cD_t d = P_d {\sf a} +\mu_V Vd +\mu_D P_d Dd$$
in the stretched case. The couple stress here is $\Pi =\rho\partial_{\nabla d}\psi$, and the entropy production now reads as
\begin{align*}\theta r &=[\alpha_0|\nabla\theta|_2^2+ \alpha_1(d|\nabla \theta)^2]/\theta + 2\mu_s|D|_2^2 + \mu_b|{\rm div}\, u|^2 \\
&+\frac{1}{\gamma}|P_d({\sf a}-\mu_PDd)|_2^2+\mu_L |P_d Dd|^2 +\mu_0 (Dd|d)^2.
\end{align*}
Summarizing, the complete model in the case of non-isotropic elasticity becomes
\begin{align}\label{cp-eqns-noniso}
\partial_t \rho +{\rm div}(\rho u)&=0\quad &\mbox{in } \Omega,\nn\\
\rho(\partial_t +u\cdot\nabla)u +\nabla \pi &= {\rm div}\, S \quad &\mbox{in } \Omega,\nn\\
\rho(\partial_t +u\cdot\nabla)\epsilon +{\rm div}\, q &= S:\nabla u -\pi{\rm div}\, u+{\rm div}(\rho\partial_{\nabla d}\psi\cD_td)\quad &\mbox{in } \Omega,\\
\gamma (\partial_t+u\cdot\nabla)d- P_d\big( {\rm div}(\rho\frac{\partial\psi}{\partial\nabla d})-\rho\nabla_d\psi\big)
&= \mu_V Vd +\mu_D P_d Dd,\quad &\mbox{in } \Omega,\nn\\
u=0,\quad q\cdot \nu=0,\quad \partial_\nu d&=0 \quad &\mbox{on } \partial\Omega.\nn\\
\rho(0)=\rho_0,\quad u(0)=u_0,\quad \theta(0)=\theta_0,\quad d(0)&=d_0\quad &\mbox{in } \Omega\nn
\end{align}
These equations have to be supplemented by the thermodynamical laws
\begin{align}\label{cp-thdyn-noniso}
\epsilon = \psi +\theta \eta, \quad \eta = -\partial_\theta \psi, \quad \kappa = \partial_\theta \epsilon,\quad \pi =\rho^2 \partial_\rho\psi,
\end{align}
and by the constitutive laws
\begin{align}\label{cp-const-noniso}
S &= S_N + S_E + S_L^{stretch} + S_L^{diss}, \nn\\
S_N&= 2\mu_s D + \mu_b {\rm div}\, u \, I, \quad S_E=  -\rho \frac{\partial \psi}{\partial\nabla d}[\nabla d]^{\sf T},\\
S_L^{stretch}& = \frac{\mu_D+\mu_V}{2\gamma}{\sf n} \otimes d + \frac{\mu_D-\mu_V}{2\gamma} d \otimes {\sf n},\quad {\sf n}= \mu_V Vd +\mu_D P_dDd-\gamma\cD_td,\nn\\
S_L^{diss}& =\frac{\mu_P}{\gamma} ({\sf n}\otimes d + d\otimes {\sf n}) + \frac{\gamma\mu_L+\mu_P^2}{2\gamma}(P_dDd\otimes d + d\otimes P_dDd) +\mu_0 (Dd|d)d\otimes d,\nn\\
q&= -\alpha_0 \nabla\theta -\alpha_1(d|\nabla\theta)d.\nn
\end{align}
Here all coefficients $\mu_j,\alpha_j$ and $\gamma$ are functions of $\rho,\theta,\tau$. For thermodynamic consistency we require as before only \eqref{cp-consist}.
We also note that the natural boundary condition at $\partial\Omega$ here becomes
$$ \nu_i\nabla_{\partial_i d}\psi = 0.$$
 Observe that this condition is fully nonlinear, in general, in contrast to the isotropic case.

Concluding, we mention as an example the classical {\em Oseen-Frank} free energy density for the isothermal incompressible case, which is given by
$$ \psi^{FO} = k_1({\rm div}\, d)^2 +k_2|d\times(\nabla\times d)|_2^2 + k_3 |d\cdot(\nabla \times d)|^2 + (k_2+k_4)[ {\rm tr}(\nabla d)^2-({\rm div}\, d)^2],$$
where $k_i$ are given constants.

\bigskip

\section{Thermodynamical Consistency and Stability}
In this section we determine the equilibria set of the complete system described above in Section 2.7, show that the critical points of the entropy functional coincide with these
equilibria and prove that they are thermodynamically stable.  We begin investigating the set of equilibria.

\medskip

\noindent
{\bf 1. Equilibria}\\
Suppose that in some time interval $t\in (t_1,t_2)$ we have $\partial_t{\sf N}(t)=0$. Then $r\geq0$ implies $r(t,x)=0$ in $\Omega$. This yields $\nabla\theta(t,x)=0$ in $\Omega$
 as $\alpha_0>0$ and $\alpha_0+\alpha_1>0$. Hence, $\theta=\theta_*$ is constant in $\Omega$.

Next, by $\mu_s>0$, $2\mu_s+n\mu_b>0$, we also have $D=0$ in $\Omega$. By Korn's inequality and the no-slip boundary condition for $u$, we hence obtain $u=u_*=0$ in $\Omega$, $t\in (t_1,t_2)$.
Therefore $\partial_t \rho=\partial_t u =0$, which implies $\nabla\pi=0$.

Finally, $\gamma>0$ yields $\cD_t d=0$ in $\Omega$, which implies that $d$
satisfies the nonlinear eigenvalue problem
\begin{align}
  \label{nlevp}
  \begin{array}{rllll}
   {\rm div}(a(x)\nabla) d +a(x)|\nabla d|_2^2 d &=0 &&\text{in } \Omega,\\
    |d|_2&=1 &&\text{in } \Omega,\\
    \partial_\nu d &=0 &&\text{on } \partial\Omega,
  \end{array}
\end{align}
where $a(x)= \lambda(\rho(t,x),\theta(t), \tau(t,x))$, for each fixed $t\in (t_1,t_2)$.
But, as the next lemma shows, this implies $\nabla d=0$ in $\Omega$. Hence, $d=d_*$ is constant.

\medskip

\noindent
{\bf Lemma 1.}\label{lma:dconst} {\em Let $q>n$, $a\in H^1_q(\Omega)$, $a>0$ and suppose that $d\in H^2_q(\Omega;\RR^n)$ satisfies \eqref{nlevp}.
Then $d$ is constant in $\Omega$.
}

\medskip

\noindent
{\em Proof.} \,
The idea is to reduce inductively the dimension $N=n$ of the vector $d$. This can be achieved by introducing polar coordinates according to
$$ d_1=c_1\cos\varphi,\; d_2=c_1\sin\varphi,\; d_j = c_{j-1},\quad j\geq3.$$
Simple computations yield
$$ 1=|d|_2^2 = |c|_2^2, \quad |\nabla d|_2^2 = |\nabla c|_2^2 + c_1^2|\nabla\varphi|_2^2, $$
and
$${\rm div}(a\nabla) c_j +a[|\nabla c|_2^2+c_1^2|\nabla\varphi|_2^2]c_j=0 \quad \mbox{ in } \Omega,$$
as well as $\partial_\nu c_j=0$ on $\partial\Omega$ for $j=2,\ldots,n-1$.
Moreover, by some more calculations we further obtain
$$-{\rm div}(a\nabla) c_1 +ac_1|\nabla \varphi|_2^2= a[|\nabla c|_2^2+c_1^2|\nabla\varphi|_2^2]c_1\quad \mbox{ in } \Omega,$$
and
$$c_1{\rm div}(a\nabla) \varphi + 2a\nabla c_1\cdot\nabla\varphi=0\quad \mbox{ in } \Omega,$$
as well as
$$ \partial_\nu c_1= c_1\partial_\nu\varphi=0\quad \mbox{ on } \partial\Omega.$$
Multiplying the second of the last equations by $c_1\varphi$ and integrating over $\Omega$ we deduce
\begin{align*}
0&=\int_\Omega [c_1{\rm div}(a\nabla)\varphi + 2a \nabla c_1\cdot\nabla \varphi)]c_1\varphi dx
= \int_\Omega {\rm div}[c_1^2a\nabla\varphi]\varphi dx = -\int_\Omega c_1^2 a|\nabla\varphi|^2dx.
\end{align*}
Hence, $c_1\nabla\varphi=0$ as $a>0$ by assumption. This implies that $c$ satisfies equation \eqref{nlevp}, where the vector $c$ has dimension $N-1$. Inductively, we arrive at dimension $N=1$ and
if $d$ is a solution of \eqref{nlevp} with dimension 1, then $d=1$ or $d=-1$ by the connectedness of $\Omega$.\hfill $\Box$

\bigskip

\noindent
Knowing that $\theta$ and $d$ are constant in $\Omega$, and $\nabla\pi=0$, we see that $\pi = \rho^2\partial_\rho\psi(\rho,\theta,0)$ is constant, hence $\rho=\rho_*$ is constant, provided
the function $\rho\mapsto \pi(\rho,\theta,0)$ is strictly increasing. This shows that we are at an equilibrium $(\rho_*,u_*,\theta_*,d_*)\in\cE$ with
$$\cE =\{(\rho_*,u_*,\theta_*,d_*)\in (0,\infty)\times\{0\}\times(0,\infty)\times \RR^n:\, |d_*|_2=1\},$$
the set of physical equilibria. In particular, the functional $-{\sf N}$ is a strict Lyapunov functional.

Observe that $\cE$ forms an $n+1$-dimensional manifold. If we take into account conservation of mass and energy,
$$ {\sf M}_0:=\int_\Omega \rho dx= \rho_*|\Omega|,\quad {\sf E}_0:= \int_\Omega(\rho|u|_2^2/2 + \rho\epsilon) dx = \rho_* \varepsilon_*|\Omega|,$$
at an equilibrium, then the values of $\rho_*$ and $\theta_*$ are uniquely determined by
$$\rho_* = {\sf M}_0/|\Omega|,\quad \epsilon_* := \epsilon(\rho_*,\theta_*,0)= {\sf E}_0/{\sf M}_0,$$
whenever $\theta\mapsto \epsilon(\rho,\theta,0)$ is strictly increasing, i.e.\ whenever  $\kappa>0$.

\bigskip

\noindent
{\bf 2. Critical Points of Total Entropy}\\
{a)} Consider the entropy functional ${\sf N}$ with constraints of prescribed mass ${\sf M} ={\sf M}_0$ and energy ${\sf E}={\sf E}_0$, as well as $G(d):=(|d|_2^2-1)/2=0$.
Suppose we have a sufficiently smooth critical point $(\rho,u,\theta, d)$ of ${\sf N}$ with $\rho,\theta>0$, subject to the constraints. Then the method of Lagrange multipliers yields $\kappa_M,\kappa_E\in\RR$ and $\kappa_G\in L_2(\Omega)$ such that
$$ \langle {\sf N}^\prime+\kappa_M {\sf M}^\prime+ \kappa_E {\sf E}^\prime +\kappa_G G^\prime|z\rangle =0,$$
where $z=(\sigma,v,\vartheta,\delta)$. We have
$$\langle {\sf M}^\prime |z\rangle =\int_\Omega \sigma dx,\quad \langle \kappa_G G^\prime|z\rangle = \int_\Omega \kappa_G d\cdot\delta dx,$$
and
$$ \langle {\sf N}^\prime|z\rangle = \int_\Omega [ (\partial_\rho(\rho\eta))\sigma + \rho\partial_\theta\eta \vartheta + \rho\partial_\tau\eta \nabla d:\nabla \delta]dx,$$
as well as
$$ \langle {\sf E}^\prime|z\rangle = \int_\Omega [ \rho u\cdot v (\partial_\rho(\rho\epsilon))\sigma + \rho\partial_\theta\epsilon \vartheta + \rho\partial_\tau\epsilon \nabla d:\nabla \delta]dx.$$
This yields the relation
\begin{align*}
0&= \int_\Omega\{ [\partial_\rho(\rho\eta) + \kappa_M + \kappa_E(\frac{1}{2} |u|_2^2 +\partial_\rho(\rho\epsilon))]\sigma + [\rho\partial_\theta\eta + \kappa_E \rho\partial_\theta\epsilon]\vartheta\} dx\\
& \quad + \int_\Omega \{\kappa_E\rho u\cdot v + [\rho\partial_\tau\eta + \kappa_E \rho \partial_\tau \epsilon]\nabla d:\nabla\delta + \kappa_G d\cdot\delta \}dx.
\end{align*}
We first vary $\vartheta$ to obtain $\rho(\partial_\theta\eta +\kappa_E\partial_\theta\epsilon)=0$, which by $\rho>0$ and by the definition of $\eta, \epsilon$ and $\kappa>0$ yields
$ \kappa_E = -1/\theta$. Hence, $\theta$ is  constant and $\kappa_E<0$. Next, varying $v$ we obtain $u=0$, as $\kappa_E$ and $\rho$ are not zero. Next we vary $\delta$, which after an integration by parts, employing the boundary condition $\partial_\nu d=0$, implies
$$ {\rm div}(\lambda\nabla) d + \kappa_G d =0\quad \mbox{ in } \Omega.$$
But then $|d|_2=1$ implies  $\kappa_G=\lambda|\nabla d|_2^2$, and $d$ is a solution of the problem \eqref{nlevp}, which by  Lemma 1 shows that $d$ is constant. Finally, we vary $\sigma$ to the result that
$ \partial_\rho(\rho\psi) = \theta \kappa_M$ is constant. As $\pi=\rho^2\partial_\rho\psi$ is strictly increasing in the variable $\rho$, this shows that $\rho$ is constant in $\Omega$ as well. Therefore, the critical points of the entropy functional are precisely the equilibria of the problem.

\bigskip

\noindent
{b)} Let
$$H:= {\sf N}^{\prime\prime}  +\kappa_E {\sf E}^{\prime\prime} $$
denote the second variation of ${\sf N}$. Note that ${\sf M}^{\prime\prime}=0$ and $\kappa_G=\lambda|\nabla d|_2^2=0$.  The identities
\begin{align*}
\rho(\partial_\tau \eta -\frac{1}{\theta}\partial_\tau \epsilon)&= -\lambda, \quad  \partial_\rho\partial_\theta(\rho\eta) -\frac{1}{\theta} \partial_\rho\partial_\theta (\rho\epsilon)=0,&\\
\partial_\rho^2(\rho\eta) -\frac{1}{\theta} \partial_\rho^2(\rho\epsilon) &= -\frac{\partial_\rho\pi}{\rho\theta},\quad  \partial_\theta^2 \eta -\frac{1}{\theta} \partial_\theta^2 \epsilon = -\frac{\kappa}{\theta^2},&
\end{align*}
imply
$$ -\langle Hz|z\rangle =\int_\Omega \big[ \frac{\partial_\rho\pi}{\rho\theta} \sigma^2 + \frac{\kappa}{\theta^2} \vartheta^2 + \lambda|\nabla \delta|_2^2\big] dx \geq0,$$
by $\kappa, \lambda,\partial_\rho\pi\geq0$. This shows that the second variation of ${\sf N}$ at an equilibrium is negative semi-definite, which means that the equilibria are {\em thermodynamically stable}.

\bigskip

\noindent
{c)} Summarizing we have the following basic result

\bigskip

\noindent
{\bf Theorem 1.}\, {\em The complete model has the following properties.

\smallskip
\noindent
{i)} Along smooth solutions total mass ${\sf M}$ and energy ${\sf E}$ are preserved.\\
{ii)} Along smooth solutions the total entropy ${\sf N}$ is non-decreasing.\\
{iii)} The negative total entropy is a strict Lyapunov functional.\\
{iv)} The condition $|d|_2=1$ is preserved along smooth solutions.\\
{v)} The equilibria are given by the set of constants
$$\cE=\{(\rho_*,0,\theta_*,d_*):\, \rho_*,\theta_*\in(0,\infty), \; d_*\in\RR^n,\; |d_*|_2=1\}.$$
Here $\rho_*, \theta_*$ are uniquely determined by  the identities
$$ \rho_*={\sf M}_0/|\Omega|,\quad  \epsilon(\rho_*,\theta_*,0) = {\sf E}_0/{\sf M}_0.$$
{vi)} The equilibria are precisely the critical points of the total entropy with prescribed mass and energy.\\
{vii)} The second variation of ${\sf N}$ with given mass and energy at equilibrium is negative semidefinite.

\smallskip\noindent
In particular, the model is {\bf thermodynamically consistent} and it is also {\bf thermodynamically stable}.}

\bigskip

\noindent
{\bf 3. The Isothermal Case}\\
In the isothermal case we set $\theta=const$ and ignore the equation for the energy. In this case, instead of the total mass specific energy $e=|u|_2^2/2 +\epsilon$, we employ the {\em available energy} $e_a$ which is defined by $e_a = |u|_2^2/2 +\psi$. We have the following balance of $e_a$ which is a direct consequence of balance of total energy and entropy
$$ \rho(\partial_t+u\cdot\nabla)e_a +{\rm div}(\Phi_e-\theta\Phi_\eta) = -\theta r -\rho\eta\cD_t\theta -\Phi_\eta\cdot\nabla\theta.$$
In the case where $\theta$ is constant this reduces to
$$ \rho(\partial_t+u\cdot\nabla)e_a +{\rm div}(\Phi_e-\theta\Phi_\eta) = - r_a,$$
with
$$r_a = 2\mu_s |D|_2^2 + \mu_b|{\rm div}\, u|^2 +\mu_0(Dd|d)^2 + \mu_L|P_dDd|_2^2 + |P_d({\rm div}(\lambda\nabla)d-\mu_PDd)|_2^2/\gamma.$$
Therefore, in the isothermal case, the total available energy ${\sf E}_a$ is a strict Ljapunov functional for the system, i.e.
$$\partial_t {\sf E}_a(t) =-\int_\Omega r_a(t,x) dx, \quad {\sf E}_a(t) = \int_\Omega \rho(t,x) e_a(t,x) dx.$$
As a consequence, the equilibrium set is the same as in the non-isothermal case, dropping temperature, hence is a manifold of dimension $n$, and when we incorporate preserved mass it is isomorphic to the unit sphere
in $\RR^n$. In this case the equations read
\begin{align}\label{isothermal}
\partial_t \rho +{\rm div}(\rho u)&=0\quad &\mbox{in } \Omega,\nn\\
\rho(\partial_t + u\cdot\nabla)u +\nabla\pi &= {\rm div} S  \quad  &\mbox{in } \Omega,\\
\gamma (\partial_t+u\cdot\nabla)d -\mu_V Vd- {\rm div}( \lambda\nabla)d&=\lambda|\nabla d|_2^2d  + \mu_D P_d Dd \quad &\mbox{in } \Omega \nn
\end{align}
where $\psi=\psi(\rho,\tau)$,  $\pi=\rho^2\partial_\rho\psi$, $\lambda = \rho\partial_\tau\psi/\theta$, and
\begin{align}\label{isoth-stress}S &= 2\mu_s D + \mu_b {\rm div}\, u \, I -\theta\lambda \nabla d[\nabla d]^{\sf T}+  S_L^{stretch}+ S_L^{diss},\nn\\
S_L^{stretch}& = \frac{\mu_D+\mu_V}{2\gamma} {\sf n} \otimes d + \frac{\mu_D-\mu_V}{2\gamma} d \otimes {\sf n},\quad {\sf n}= \mu_V Vd +\mu_D P_dDd-\gamma\cD_t d.\\
S_L^{diss}& =\frac{\mu_P}{\gamma}({\sf n}\otimes d + d\otimes {\sf n}) + \frac{\gamma\mu_L+\mu_P^2}{2\gamma}(P_dDd\otimes d + d\otimes P_dDd) +\mu_0 (Dd|d)d\otimes d,\nn\\
 \end{align}
If one further restricts to the incompressible case $\rho=const >0$, $\lambda$ constant, $\mu_0=\mu_D=\mu_V=\mu_P=\mu_L=0$, with $\mu=\mu_s$ one obtains
the so-called {\em isothermal simplified Ericksen-Leslie model}
\begin{align}\label{simplified}
\rho(\partial_t + u\cdot\nabla)u +\nabla\pi& = \mu_s \Delta u -\lambda {\rm div}( \nabla d[\nabla d]^{\sf T})  \quad &\mbox{in } \Omega,\nn\\
|d|_2 = 1, \quad
{\rm div}\, u &=0 \quad &\mbox{in }\Omega,\\
\gamma (\partial_t+u\cdot\nabla)d -\lambda \Delta d&=\lambda|\nabla d|^2d \quad &\mbox{in }\Omega.\nn
\end{align}
Of course, in all cases we have to add initial conditions as well as  boundary conditions $u=\partial_\nu d=0$ on $\partial\Omega$.  Problem \eqref{simplified} subject to the condition
$|d|_2=1$ in $\Omega$ has been analyzed in a fairly complete manner  in the recent article \cite{HNPS14} by Hieber, Nesensohn, Pr\"uss and  Schade.

\section{Analysis of the Non-Isothermal Simplified Model}

In this section we consider the {\em incompressible case} $\rho=const$ and we let $\mu=\mu_s$. Hence, the pressure $\pi$ is no longer determined by Maxwell's relation;
it is now a free variable, a Lagrangian multiplier to cover the constraint ${\rm div}\, u=0$.
Furthermore, in the following  we neglect stretching, i.e.\ we assume  $\mu_D=\mu_V=\mu_P=\mu_L=\mu_0=0$.  For simplicity we also set $\alpha_1=0$ and $\alpha=\alpha_0$.
Then the resulting model - which we call the {\em non-isothermal simplified Ericksen-Leslie model} - reads as follows.
\begin{align}\label{reduced-model}
\rho\cD_t u -2{\rm div}(\mu D) +\nabla \pi &= {\rm div}(\lambda\nabla d[\nabla d]^{\sf T})\quad &\mbox{in } \Omega,\nn\\
|d|_2 =1,\quad {\rm div}\; u &=0\quad &\mbox{in } \Omega,\nn\\
\rho\kappa \cD_t \theta -{\rm div}(\alpha\nabla\theta)&= 2\mu|D|_2^2 +\theta\lambda\nabla d[\nabla d]^{\sf T}:D\quad &\nn\\
& \quad - \rho\partial_\tau\epsilon \nabla d:\cD_t\nabla d +{\rm div}(\lambda\nabla d \cD_t)\quad &\mbox{in } \Omega,\\
\gamma \cD_t d -{\rm div}(\lambda\nabla)d &= \lambda |\nabla d|^2_2 d\quad &\mbox{in } \Omega,\nn\\
u=\partial_\nu\theta=\partial_\nu d&=0\quad &\mbox{on } \partial\Omega,\nn\\
u(0)=u_0,\; \theta(0)=\theta_0,\; d(0)&=d_0\quad &\mbox{in } \Omega.\nn
\end{align}
Recall that $\rho>0$ is constant and $\alpha,\gamma,\mu$ as well as $\lambda=\rho\partial_\tau\psi/\theta, \kappa =\partial_\theta\epsilon=-\theta \partial_\theta^2\psi$ are functions of $\theta>0$ and $\tau\geq0$.

\bigskip

\noindent
{\bf 1. Regularity Assumptions {\bf (R)}}\\
The parameter functions should have the following minimal regularity properties:
$$ \mu, \alpha, \gamma\in C^2((0,\infty)\times[0,\infty)),\quad \psi\in C^4((0,\infty)\times [0,\infty));$$
We also require the positivity conditions
$$  \mu>0,\quad \alpha>0, \quad \kappa>0, \quad \gamma>0,\quad \lambda>0,$$
which have been mentioned before, but
for well-posedness of the problem for $d$ we need in addition to require
$ \lambda +2\tau \partial_\tau \lambda>0$. We assume that $\Omega\subset\RR^n$ is a bounded domain with $C^{3-}$-boundary.

\bigskip

\noindent
{\bf 2. Maximal $L_p$-Regularity of the Principal Linearization}\\
 The equation for $u$ will turn out to be only weakly coupled, so we first concentrate on the system for $w:=[\theta,d]^{\sf T}$. The {\em principal part of the linearization} becomes
\begin{eqnarray}
\partial_t w +\cA(w_0,\nabla) w &=&f\qquad \mbox{ {in} } \Omega,\nn\\
\partial_\nu w &=&0 \qquad \mbox{ {on} } \partial\Omega,\\
w(0)&=&w_0\quad\,\, \mbox{ in } \Omega.\nn
\end{eqnarray}
The matrix $\cA=\cA(w_0,\nabla)$ reads as
$$ \cA= \left[
\begin{array} {cc}
-a_0 {\Delta} - a_1\nabla d_0[\nabla d_0]^{\sf T}:{\nabla^2}, & b_0 \nabla d_0: ({\lambda_0\Delta}  + \partial_\tau\lambda_0[\nabla d_0]^{\sf T}\nabla d_0:{\nabla^2})\nabla\\
b_1[\nabla d_0]^{\sf T} {\nabla}, & -\gamma_0^{-1} (\lambda_0 {\Delta}  + \partial_\tau\lambda_0 [\nabla d_0]^{\sf T}\otimes\nabla d_0:{\nabla^2}) .
\end{array}
\right].$$
Here  $\kappa_0=\kappa(\theta_0,\tau_0)$ etc., and we used the abbreviations
$$ a_0= \frac{\alpha_0}{\rho\kappa_0},\quad a_1= \frac{\rho\theta_0[\partial_\tau\eta_0]^2}{\gamma_0\kappa_0},\quad b_0 =
\frac{\theta\partial_\tau \eta_0}{\gamma_0\kappa_0},\quad b_1=\frac{\rho\partial_\tau\eta_0}{\gamma_0}.$$
Note that {$\cA(w_0,\nabla)$} is second order in the diagonal, but third and first order off-diagonal! This is a mixed-order problem subject to Neumann boundary conditions
and subject to variable, non-smooth coefficients. For resolvent estimates within the $L^p$-setting for various mixed-order systems we refer to the work of Grubb \cite{Gru77}.
Regarding the {\em maximal $L^p$-regularity} of this nonstandard problem, we do not know of any general theory covering the above situation. However, we note for the whole space case there is the theory of Denk and Kaip \cite{DeKa13} available.
 Nevertheless, we prove
{maximal $L_p$-regularity} for this problem in the following.

To this end, fix $q\in(1,\infty)$ and choose as a base space
$$Y_0:= L_q(\Omega)\times H^1_q(\Omega;\RR^n),$$
and as a regularity space
$$ Y_1:=\{w=(\theta,d)\in H^2_q(\Omega)\times H^3_q(\Omega;\RR^n):\; \partial_\nu \theta=\partial_\nu d=0 \mbox{ {on} } \partial\Omega\},$$
equipped with their natural norms.
We will also employ the {\em time-weighted spaces} defined by
$$ {\sf y} \in H^m_{p,\mu}(J;Y) \quad \Leftrightarrow \quad t^{1-\mu}{\sf y} \in H^m_p(J;Y),\quad m\in{\mathbb N}_0,\; \mu\in (1/p,1].$$
The trace space $Y_{\gamma,\mu}$ is then given by
$$Y_{\gamma,\mu}=\{ (\theta,d)\in B_{qp}^{2(\mu-1/p)}(\Omega)\times B_{qp}^{2(\mu-1/p)+1}(\Omega;\RR^n):\; \partial_\nu d=0 \mbox{ on } \partial\Omega\},$$
provided
$$ \frac{1}{p}<\mu < \frac{1}{2}+\frac{1}{p}+\frac{1}{2q};$$
otherwise one has to add $\partial_\nu\theta=0$ in the definition of $Y_{\gamma,\mu}$. In order  to  profit from the embedding
\begin{equation}\label{Y-embedding}
Y_{\gamma,\mu} \hookrightarrow C(\overline{\Omega}) \times  C^1(\overline{\Omega};\RR^n),
\end{equation}
we will always assume
$$ 1\geq \mu >\frac{1}{p}+\frac{n}{2q}.$$
Then by means of the assumptions stated before we obtain the following result.

\bigskip

\noindent
{\bf Theorem 2.} {\em Assume {(R)}, $1/p+n/2q<\mu\leq 1$, and suppose that {$w_0\in Y_{\gamma,\mu}$}. Then the differential operator {$A_2(w_0)$} defined by
{$A_2(w_0)w:=\cA(w_0,\nabla)w$} with domain {${\sf D}(A_2(w_0)):=Y_1$} has maximal {$L_p$}-regularity in {$Y_0$} and thus  also
maximal $L_{p,\mu}$-regularity in $Y_0$. }

\bigskip

\begin{proof}
The proof is based on the results and techniques developed by  Denk, Hieber, Pr\"{u}ss in \cite{DHP03} for the case $\mu=1$.
By the results due to Pr\"{u}ss and Simonett \cite{PS04},  these results extend to general $\mu\in (1/p+n/2q,1]$, as the coefficients have enough regularity by the
embedding \eqref{Y-embedding}.

\medskip

\noindent
{a)} The case $\Omega=\RR^n$ with  constant coefficients\\
In the sequel, we denote the covariable for $t$ by $z$ and that for $x$ by $\xi$.
The symbol $\cA(\xi)$ of $\cA(w_0,\nabla)$ reads as
\begin{align*}
\cA(\xi) = \left[ \begin{array}{cc}
a_0|\xi|^2 + a_1 |c(\xi)|^2 & -i b_0(\lambda_0|\xi|^2+\partial_\tau\lambda_0 |c(\xi)|^2)c(\xi)^{\sf T}\\
i b_1c(\xi)& \frac{\lambda_0}{\gamma_0} |\xi|^2 +\frac{\partial_\tau\lambda_0}{\gamma_0}c(\xi)\otimes c(\xi)
\end{array}\right],
\end{align*}
where $c(\xi) = \xi\cdot\nabla d_0$. It is convenient to reduce this symbol for the variable $w_{red}=[\theta, d_{red}]^{\sf T}$ where $d_{red}=c(\xi)\cdot d$. The reduced symbol $\cA_{red}(\xi)$ becomes
\begin{align*}
\cA_{red}(\xi) = \left[ \begin{array}{cc}
a_0|\xi|^2 + a_1 |c(\xi)|^2 & -i b_0(\lambda_0|\xi|^2+\partial_\tau\lambda_0 |c(\xi)|^2)\\
i b_1|c(\xi)|^2& \frac{\lambda_0}{\gamma_0} |\xi|^2 +\frac{\partial_\tau\lambda_0}{\gamma_0}|c(\xi)|^2
\end{array}\right].
\end{align*}
This symbol is homogeneous of second order and  {\em not} strongly elliptic. However, it is  {\em normally elliptic} in the sense of \cite{DHP03}  as its spectrum
satisfies $\sigma(\cA_{red}(\xi))\subset (0,\infty)$ for each $\xi\neq0$. The latter can be seen
by considering
\begin{align*}
{\rm det}(z+\cA_{red}(\xi))&= (z+a_0|\xi|^2 +a_1|c(\xi)|^2)( z +  \frac{\lambda_0}{\gamma_0} |\xi|^2 +\frac{\partial_\tau\lambda_0}{\gamma_0}|c(\xi)|^2)\\
 &\quad -b_0b_1|c(\xi)|^2(\frac{\lambda_0}{\gamma_0} |\xi|^2 +\frac{\partial_\tau\lambda_0}{\gamma_0}|c(\xi)|^2),\quad \xi\neq0,
 \end{align*}
which has two negative zeros, as $\alpha,\gamma,\kappa>0$ and $\lambda+ 2\tau\partial_\tau\lambda>0$.  Therefore, by Section 6 of \cite{DHP03}, the $L_p$-realization
$A_{red}$ of $\cA_{red}$ has maximal $L_p$-regularity. This shows that whenever $f_\theta\in L_p(J;L_q(\RR^n))$ and $f_d\in L_p(J;H^1_q(\Omega;\RR^n))$ are given, there is a unique solution $$w_{red}=[\theta,\bar{d}]^{\sf T}\in {_0H}^1_p(J;L_q(\RR^n;\RR^2))\cap L_p(J;H^2_q(\RR^n;\RR^2))$$ of
$$ \partial_tw_{red}+A_{red}w_{red} = f_{red},\quad t>0,\; w_{red}(0)=0,$$
where $f_{red}=[f_\theta,-ic(\nabla)\cdot f_d]^{\sf T}$. To obtain $d$, it remains to solve the problem
$$ \partial_t d -\frac{\lambda_0}{\gamma_0} \Delta d = f^1_d:= f_d +i \frac{\partial_\tau\lambda_0}{\gamma_0} c(\nabla) d_{red}- b_1 c(\nabla)\theta,\quad t>0,\; d(0)=0,$$
with maximal $L_p$-regularity of $-\Delta$ to obtain a unique solution
$$d\in {_0H}^1_p(J;H^1_q(\RR^n;\RR^n))\cap L_p(J;H^3_q(\RR^n;\RR^n)),$$
as $\lambda,\gamma>0$ and $f^1_d\in L_p(J;H^1_q(\RR^n;\RR^n))$. This proves Theorem 2 in the case $\Omega=\RR^n$ with constant coefficients. As detailed in Section 6 of
\cite{DHP03}, this assertion extends by perturbation and localization to variable coefficients, still in the case $\Omega=\RR^n$.

\medskip

\noindent
{b)} The case $\Omega=\RR^n_+$ with  constant coefficients\\
It is convenient to replace $x\in\RR^n_+$ by $(x,y)\in \RR^{n-1}\times\RR_+$.
On the symbolic level we have to replace $\xi$ by $\xi-i\nu\partial_y$, where $\nu$ denotes the outer normal at a boundary point of $\Omega$, and $\xi\cdot\nu=0$. Then $c(\xi)$ becomes
$$ c(\xi-i\nu\partial_y) = (\xi-i\nu\partial_y)\cdot \nabla d_0 = \xi\cdot\nabla d_0 -i \nu\cdot\nabla d_0\partial_y = \xi\cdot\nabla d_0,$$
as $\partial_\nu d_0=0$. Therefore the symbol $\cA(\xi)$ from Step {a)} is replaced by
$$ -E\partial_y^2 + z+\cA(\xi),$$
where
$$ E =\left[\begin{array}{cc}
a_0& -ib_0 c(\xi)^{\sf T}\\
0& (\lambda_0\gamma_0)I
\end{array}\right].$$
Considering again the reduced variables $w_{red}=(\theta, d_{red})$ with $d_{red}=c(\xi)\cdot d$, the reduced symbol becomes
$$ -E_{red}\partial_y^2 + z+\cA_{red}(\xi),$$
with
$$ E =\left[\begin{array}{cc}
a_0& -ib_0\\
0& \lambda_0\gamma_0
\end{array}\right].$$
To apply the half-space theory for normally elliptic operators in Section 7 of  \cite{DHP03}, we need  to verify the
corresponding Lopatinskii-Shapiro condition {(LS)} which states the following.

\medskip

\noindent
{\em If $w\in C_0(\RR_+;\RR^2)$ satisfies
$$ \text{(LS)}\quad -E_{red}\partial_y^2w(y) + (z+\cA_{red}(\xi))w(y)=0,\quad y>0, \quad \partial_y w(0)=0,$$
then $w=0$.}

\medskip

\noindent
Condition (LS) can be proved without much pain. In fact, we observe that $E_{red}^{-1}(z+\cA_{red}(\xi))$ has no eigenvalues $-\omega^2\leq 0$; otherwise $z\in\CC\setminus(-\infty,0]$ would be a solution of
$$ 0 = {\rm det}(z+\cA_{red}(\xi)+\omega^2E_{red}) = {\rm det}(z+\cA_{red}(\xi+\omega\nu)),$$
which by Step {a)} is impossible.  Therefore $B=(E_{red}^{-1}(z+\cA_{red}(\xi)))^{1/2}$ is well-defined and has spectrum in $\CC_+$. Thus, $w_{red}(y):= e^{-By}w_b$ is the unique stable solution of
$$-E_{red}\partial_y^2w_{red}(y) + (z+\cA_{red}(\xi))w_{red}(y)=0,\quad y>0, \quad  w_{red}(0)=w_b.$$
The Neumann condition implies
$$0 =\partial_y w_{red}(0) = -B w_b,$$
hence $w_b=0$, as $B$ is invertible, for all $(z,\xi)\neq(0,0)$, $z\in \CC\setminus(-\infty,0]$, $\xi\in \RR^{n}$, $\xi\cdot\nu=0$. Therefore,
the techniques of Section 7  in \cite{DHP03} apply and  show  that the reduced problem has maximal $L_p$-regularity in $\RR^n_+$.

As a result, given $f_\theta\in L_p(J;L_q(\RR^n_+))$, $f_d\in L_p(J;H^1_q(\RR^n_+;\RR^n))$, with $f_{red}=[f_\theta,-ic(\nabla_x)\cdot f_d]^{\sf T}$,
we find a unique solution $w_{red}=[\theta,d_{red}]^{\sf T}$ of the problem
\begin{align*}
\partial_t w_{red} + A_{red} w_{red} &= f_{red} \quad \mbox{ in } \RR^n_+,\quad w_{red}(0)=0,
\end{align*}
within the class
$$ w_{red}\in {_0H}^1_p(J;L_q(\RR^n_+;\RR^2)) \cap L_p(J;H^2_q(\RR^n_+;\RR^2)),$$
where $A_{red}$ denotes the realization of $\cA_{red}$ in $L_q(\RR^n_+;\RR^2)$ with Neumann boundary condition.

Next, we solve the remaining problem for $d$
\begin{align*}
\partial_t d -\frac{\lambda_0}{\gamma_0} \Delta  d &= f^1_d \quad \mbox{ in } \RR^n_+,\; d(0)=0,\\
\partial_\nu d &=0 \quad \mbox{ on } \partial\RR^n_+
\end{align*}
with
$$ f^1_d := f_d +i \frac{\partial_\tau\lambda_0}{\gamma_0} c(\nabla) d_{red}- b_1 c(\nabla)\theta ,$$
in a similar way as in Step {a)} by employing maximal $L_p$-regularity for $-\Delta$. This yields a unique solution
$$d\in {_0H}^1_p(J;L_q(\RR^n_+;\RR^n)) \cap L_p(J;H^2_q(\RR^n_+;\RR^n)).$$
As $f^1_d \in L_p(J;H^1_q(\RR^n_+;\RR^n))$ we may differentiate the equation for $d$ tangentially to obtain also
$$\nabla_xd\in {_0H}^1_p(J;L_q(\RR^n_+;\RR^{(n-1)\times n})) \cap L_p(J;H^2_q(\RR^n_+;\RR^{(n-1)\times n})).$$
On the other hand, we may also take the derivative with respect to the normal variable $y$ in order  to obtain a  problem with Dirichlet boundary conditions for $v:=\partial_y d$.
We solve this with maximal $L_p$-regularity to obtain
$$\partial_yd\in {_0H}^1_p(J;L_q(\RR^n_+;\RR^n)) \cap L_p(J;H^2_q(\RR^n_+;\RR^n)).$$
This proves Theorem 2 for the case $\Omega=\RR^n_+$ with constant coefficients. As described in \cite{DHP03}, Section 7, this assertion extends by perturbation and
localization to variable coefficients, still in the case $\Omega=\RR^n_+$.

\medskip

\noindent
{c)} General domains and variable coefficients\\
Here we follow the line given in \cite{DHP03}, Section 8.  We may use a perturbation argument to extend the result for the half-space to a bent half-space and then employ
the localization method to prove Theorem 2 for general domains with $C^{3-}$-boundary.
\end{proof}

\bigskip

\noindent
{\bf 3. Local-Wellposedess}\\
We rewrite the above problem as an abstract quasilinear evolution equation of the form
\begin{equation} \label{NLCF}
\dot{z} +A(z) z = F(z),\quad t>0,\; z(0)=z_0.
\end{equation}
Here $z=(u,w)=(u,\theta,d)$ and we apply the {\em Helmholtz projection} $\PP$ to the equation for $u$.  The base space will be $X_0:= L_{q,\sigma}(\Omega)\times Y_0$, where the subscript $\sigma$ means {\em solenoidal}. Then with the {\em generalized Stokes operator} $A_1(w)=-\PP\mu(\theta,\tau)\Delta$, we define the regularity space by
$$ X_1:={\sf D}(A_1)\times Y_1,\quad {\sf D}(A_1) =\{ u\in H^2_q(\Omega;\RR^n)\cap L_{q,\sigma}(\Omega):\, u=0 \mbox{ {on} } \partial\Omega\}.$$
The operator  $A(z)$ is defined by
$ A(z) ={\rm diag}(A_1(w),A_2(w)),$
and $F(z)$ collects all lower order terms.

In order to prove local well-posedness of (\ref{NLCF}), we may now resort to abstract theory, e.g.\ to the results by  K\"ohne, Pr\"uss, Wilke in \cite{KPW10} and by
LeCrone, Pr\"{u}ss and Wilke in \cite{LPW14}.

Then, by Theorem 2 and by the maximal regularity of the generalized Stokes operator, see e.g.\ Bothe and Pr\"uss \cite{BP07}, $A(z)$ has maximal $L_p$-regularity.
For an interval $J=[0,a]$, the {\em solution space} $\EE_\mu(J)$ will be
$$ \EE_\mu(J)= H^1_{p,\mu}(J;X_0)\cap L_{p,\mu}(J;X_1).$$
The time-trace space
$X_{\gamma,\mu}$ of $\EE_\mu(J)$ is given by
$$X_{\gamma,\mu} = \{ u\in B_{qp}^{2(\mu-1/p)}(\Omega)^n\cap L_{q,\sigma}(\Omega):\, u_{|_{\partial\Omega}}=0\}\times Y_{\gamma,\mu};$$
it satisfies
$$ X_{\gamma,\mu}\hookrightarrow B_{qp}^{2(\mu-1/p)}(\Omega)^{n+1}\times B_{qp}^{1+2(\mu-1/p)}(\Omega)^n\hookrightarrow C(\overline{\Omega})^{n+1}\times C^1(\overline{\Omega})^n,$$
provided
\begin{equation}\label{pq}
\frac{1}{p} +\frac{n}{2q}<\mu\leq 1.
\end{equation}
Here $B^s_{pq}$ denote as usual the Besov spaces; see e.g.\ Triebel \cite{Tri92}.
Then $A,F$ satisfy the requirements in the paper by LeCrone, Pr\"{u}ss and Wilke \cite{LPW14}, and so we have local well-posedness. If $\frac{1}{p} +\frac{n}{2q}+\frac{1}{2}<\mu\leq 1$, the conditions of K\"{o}hne, Pr\"{u}ss, and Wilke \cite{KPW10} also hold. In particular, defining the {\em state manifold} of (\ref{NLCF}) by
$$\cSM=\{ (u,\theta, d)\in X_{\gamma}: \, \theta>0,\; |d|_2=1\},\quad X_\gamma:=X_{\gamma,1},$$
then $\cSM$ is {\em locally positive invariant} for the semi-flow, total energy ${\sf E}$ is preserved, and the negative total entropy $-{\sf N}$ is a strict Lyapunov functional for the semi-flow on $\cSM$. Summarizing we have the following result

\medskip

\noindent
{\bf Theorem 3.} \label{local}
{\em Assume {(R)}, let $p,q,\mu$ be subject to \eqref{pq}, and let $z_0 \in X_{\gamma,\mu}$. Then for some $a=a(z_0)>0$, there is a unique solution
$$z\in H^1_{p,\mu}(J,X_0)\cap L_{p,\mu}(J;X_1),\quad J=[0,a],$$
of \eqref{NLCF}, i.e.\ \eqref{reduced-model} on $J$. Moreover,
$$z\in C([0,a];X_{\gamma,\mu})\cap C((0,a];X_\gamma),$$
i.e. the solution regularizes instantly in time.
It depends continuously on $z_0$ and exists on a maximal time interval $J(z_0) = [0,t^+(z_0))$.
Moreover,
\begin{align*}t [\frac{d}{dt}] z \in  H^1_{p,\mu}(J;X_0)\cap L_{p,\mu}(J;X_1),\end{align*}
and $|d(t,x)|_2\equiv 1$, ${\sf E}(t)\equiv {\sf E}_0$, and $-N$ is a strict Lyapunov functional.
Furthermore, the problem \eqref{NLCF} generates a local semi-flow in its natural state manifold $\cSM$.
}

\bigskip

\noindent
{\bf 4. The Generalized Principle of Linearized Stability}\\
Consider the autonomous
quasilinear problem
\begin{equation}
\label{u-equation}
\dot{z}(t)+A(z(t))z(t)=F(z(t)),\quad t>0, \quad z(0)=z_0.
\end{equation}
Here we assume
\begin{equation}
\label{AF}
(A,F)\in C^1(V,\cB(X_1,X_0)\times X_0),
\end{equation}
where $V\subset X_\gamma$ is open.
Let $ \cE\subset V\cap X_1$ denote the set of  equilibrium solutions of (\ref{u-equation}), which means that
$$
z\in\cE \quad {\mbox{ if and only if }}\quad z\in V\cap X_1,
\; A(z)z=F(z).
$$
Given an  element $z_*\in\cE$,  we assume that $z_*$ is
contained in an $m$-dimensional manifold of equilibria. This means that there
is an open subset $U\subset\R^m$, $0\in U$, and a $C^1$-function
$\psi:U\rightarrow X_1$,  such that
\begin{equation}
\psi(U)\subset \cE,\quad\psi(0)=z_*,\quad
A(\psi(\zeta))\psi(\zeta)=F(\psi(\zeta)),\quad \zeta\in U.
\end{equation}
and the rank of $\psi^\prime(0)$ equals $m$.

Let $A_0$ denote the linearization of $A(z)z-F(z)$ at $z_*$, i.e.
$$ A_0h = A(z_*)h + [A^\prime(z_*)h]z_*-F^\prime(z_*)h.$$
We call $z_*\in\cE$ {\em normally stable} if the following conditions hold.

\medskip

\noindent
{(i)}\, near $z_*$ the set  $\cE$ is a $C^1$-manifold in $X_1$,  ${\rm dim}\, \cE=m\in\N_0$,\\
{(ii)} \, the tangent space for $\cE$ at $z_*$ is isomorphic to ${\sf N}(A_0)$,\\
{(iii)} \, $0$ is a semi-simple eigenvalue of
$A_0$, i.e.\ $ {\sf N}(A_0)\oplus {\sf R}(A_0)=X_0$,\\
{(iv)} \, $\sigma(A_0)\setminus\{0\}\subset \C_+=\{\zeta\in\C:\, {\rm Re}\, \zeta>0\}$.

\medskip

\noindent
The following result is due to Pr\"uss, Simonett and Zacher \cite{PSZ09}.

\bigskip

\noindent
{\bf Theorem 4.} \label{th:1}
{\em Let $1<p<\infty$. Suppose $z_*\in V\cap X_1$ is an equilibrium of (\ref{u-equation}) and that
$(A,F)$ satisfy (\ref{AF}) and that $A(z_*)$ has the property of maximal $L_p$-regularity.
Assume further that $z_*$ is normally stable.

Then $z_*$ is stable in $X_\gamma$, and there exists $\delta>0$ such
that the unique solution $z(t)$ of (\ref{u-equation}) with initial
value $z_0\in X_\gamma$ satisfying $|z_0-z_*|_{\gamma}<\delta$
exists on $\R_+$ and converges at an exponential rate in $X_\gamma$
to some $z_\infty\in\cE$ as $t\rightarrow\infty$.}

\bigskip

\noindent
It is worthwhile to note that in case $m=0$, $z_*$ is necessarily isolated by (i). Then Theorem 4 reduces to the usual {\em principal of linearized stability}, as
(ii), (iii), (iv) are equivalent to $\sigma(A_0)\subset \CC_+$.

\bigskip

\noindent
{\bf 5. Linear Stability of Equilibria}\\
The linearization of (\ref{NLCF}), i.e.\ of \eqref{reduced-model} at an equilibrium $z_*=(0,\theta_*,d_*)$ is given by the operator
$$ A_*=-{\rm diag}\big((\mu_*/\rho)\PP\Delta, (\alpha_*/\rho\kappa_*)\Delta, (\lambda_*/\gamma_*)\Delta\big)$$
in the base space $X_0$ with domain ${\sf D}(A_*)=X_1$. This operator has maximal $L_p$-regularity, it is the negative generator of a { compact analytic $C_0$-semigroup},
and it has { compact resolvent}. So its spectrum consists only of countably many eigenvalues of finite multiplicity, which are all positive, hence stable,
except for $0$. The eigenvalue $0$ is semi-simple, its eigenspace is given by
$$ {\sf N}(A_*)=\{(0,\vartheta, {\sf d}):\, \vartheta\in\RR, {\sf d}\in\RR^n\},$$
hence it coincides with the set of constant equilibria $\bar{\cE}$, when ignoring the constraint $|d|_2=1$ and conservation of energy.
Therefore each such equilibrium is normally stable.

\bigskip

\noindent
{\bf 6. Nonlinear Stability}\\
We have {\em stability with  asymptotic phase} for the equilibria of (\ref{NLCF}).

\medskip

\noindent
{\bf Theorem 5.} {\em Assume { (R)}. Then any equilibrium $z_*\in \bar{\cE}$ of (\ref{NLCF}) is stable in $X_\gamma$. Moreover, for each $z_*\in \bar{\cE}$ there is
$\varepsilon > 0$ such that if
$|z_0 -z_*|_{X_{\gamma,\mu}} \leq \ve$, then the solution $z$ of (\ref{NLCF}) with initial value $z_0$ exists globally in time and converges at an exponential rate in $X_\gamma$ to some $z_\infty\in \bar{\cE}$.
}
\medskip

\noindent
This result is proved by means of the {\em generalized principle of linearized stability}, Theorem 4,   above.
In fact, by the previous section we know that each equilibrium $z_*=(0,\theta_*,d_*)$ is normally stable.

\bigskip

\noindent
{\bf 7. Long-Time Behaviour}\\
We conclude this paper with a result on the convergence of
solutions to equilibria in the topology of the state manifold $\cSM$.

\bigskip

\noindent
{\bf Theorem 6.}
{\em  Assume  { (R)} and let  $z$ be the solution of (\ref{NLCF}), i.e.\ of \eqref{reduced-model}, with initial value  $z_0\in \cSM$. Then the following assertions hold.

\smallskip

\noindent
a) If we suppose
$$
\sup_{t \in (0,t^+(z_0))} [|z(t)|_{X_{\gamma,\mu}} + |1/\theta(t)|_{L_\infty}]<\infty,
$$
then $t^+(z_0)=\infty$ and $z$ is a global solution.\\
b) If $z$ is a global solution, bounded in $X_{\gamma,\mu}$ and with $1/\theta$ bounded, then $z$ converges exponentially in $\cSM$ to an
equilibrium $z_\infty\in \mathcal E$ of (\ref{NLCF}), as $t\to \infty$.
}

\bigskip

\noindent
This result follows from abstract dynamical system arguments involving the strict Lyapunov functional $-{\sf N}$, as well as the nonlinear stability result;
see K\"ohne, Pr\"uss,  Wilke \cite{KPW10}. Note that, by a compactness argument, the converse of b) is also valid.


\end{document}